\documentclass[preprint,authoryear,3p]{elsarticle}
\usepackage{amsmath,amsfonts,amsthm,amssymb,mathtools}
\usepackage{graphicx}
\newcommand*{\eps}{{\varepsilon}}
\newcommand*{\rd}{{\mathrm{d}}}

\newcommand*{\zb}{{\boldsymbol{z}}}
\newtheorem{theorem}{Theorem}[section]
\newtheorem{lemma}[theorem]{Lemma}
\theoremstyle{definition}

\theoremstyle{remark}
\newtheorem{remark}[theorem]{Remark}
\journal{Statistics and Probability Letters}

\begin{document}
\begin{frontmatter}
\title{Efficient almost-exact L\'evy area sampling}
\author{Simon J.A. Malham\corref{cor1}}
\ead{S.J.Malham@hw.ac.uk}
\author{Anke Wiese}
\ead{A.Wiese@hw.ac.uk}
\address{Maxwell Institute for Mathematical Sciences 
              and School of Mathematical and Computer Sciences,
              Heriot-Watt University, Edinburgh EH14 4AS, UK (+44 131 4513254)}
\cortext[cor1]{Corresponding author}
\begin{abstract}
We present a new method for sampling the L\'evy area for 
a two-dimensional Wiener process conditioned on its endpoints.
An efficient sampler for the L\'evy area is required to
implement a strong Milstein numerical scheme
to approximate the solution of a stochastic differential 
equation driven by a two-dimensional Wiener process
whose diffusion vector fields do not commute. Our method
is simple and complementary to those 
of Gaines--Lyons and Wiktorsson, and amenable to 
quasi-Monte Carlo implementation.
It is based on representing the L\'evy area
by an infinite weighted sum of independent Logistic random variables.
We use Chebychev polynomials to approximate the inverse
distribution function of sums of independent Logistic random 
variables in three characteristic regimes. The error is
controlled by the degree of the polynomials, we set
the error to be uniformly $10^{-12}$. 
We thus establish a strong almost-exact 
L\'evy area sampling method. The complexity of our method is 
square logarithmic.
We indicate how it can contribute to efficient
sampling in higher dimensions.
\end{abstract}
\end{frontmatter}

\section{Introduction}\label{sec:intro}
We consider the problem of sampling the L\'evy area for 
a two-dimensional Wiener process $(W^1_t,W^2_t)^{\text{\tiny{T}}}$ 
conditioned on its endpoints.
Indeed, on each computational timestep of size $h$, we must 
generate two independent  
sample Wiener increments, $\Delta W^1$ and $\Delta W^2$, 
and a sample of the L\'evy area 
\begin{equation*}
A(h)\coloneqq\frac{1}{2}\int_t^{t+h}\int_t^{\tau_1}
\,\rd W^1_{\tau_2}\rd W^2_{\tau_1}-\rd W^2_{\tau_2}\rd W^1_{\tau_1}.
\end{equation*}
\cite{W} proposed approximating the L\'evy area,  
given $\Delta W^1$ and $\Delta W^2$, by \citep[see][]{Levy}
\begin{equation*}
\frac{h}{2\pi}\Biggl[\sum_{n=1}^N\frac{1}{n}
\Bigl(U_{n}\bigl(Y_{n}-\sqrt{\tfrac{2}{h}}\Delta W^2\bigr)
-V_{n}\bigl(X_{n}-\sqrt{\tfrac{2}{h}}\Delta W^1\bigr)\Bigr)
+\biggl(2(1+a^2)\Bigl(\frac{\pi^2}{6}
-\sum_{n=1}^N\frac{1}{n^2}\Bigr)\biggr)^{\frac12}\,Z\Biggr].
\end{equation*}
Here $a^2\coloneqq\bigl((\Delta W^1)^2+(\Delta W^2)^2\bigr)/h$
and $U_{n},V_{n},X_{n},Y_{n}$, for $n=1,\ldots,N$, and $Z$ 
are independent standard Normal random variables. 
Without the tail term involving $Z$, 
it is the Kloeden--Platen--Wright approximation \citep{KPW}
with mean-square error of order $h^2/N$.
With Wiktorsson's tail approximation,
the mean-square error improves to $h^2/N^2$.
This method is not restricted to a
two-dimensional Wiener process.

In general, to implement a strong Milstein method
to approximate the solution of a stochastic differential
equation driven by a two-dimensional Wiener process,
we must sufficiently accurately strongly sample the 
L\'evy area. We measure the \emph{complexity} (computational effort) 
associated with such an approximation by the number of 
uniform random variables required to generate a L\'evy area 
sample on each timestep of mean-square accuracy $\epsilon$.
The smaller the complexity, the more effective is the simulation method. 
Roughly, for the Kloeden--Platen--Wright method, to achieve
accuracy $\epsilon$, we require $N$ of order $h^2/\epsilon$.
The number of uniform random variables required to 
generate an approximation truncated at $N$ terms is of order $N$.
Hence the complexity is $h^2/\epsilon$.
For Wiktorsson's method, the number of uniform samples we require 
is also of order $N$. However to achieve accuracy $\epsilon$ we require
$N$ to be of order $h/\epsilon^{1/2}$, which is a significantly improved 
complexity. The method proposed by \citet{RW} also has complexity $h/\eps^{1/2}$; 
though see Section~\ref{sec:conclu}. 
The \citet{GL1994} method is an exact acceptance-rejection method. 
Hence it cannot be used for quasi-Monte Carlo simulations. 
It is reported to be ``fast but complicated to implement'', see \citet{RW}.
\citet{SH} derive a new series representation of the joint 
distribution function. 
However in practice, a large number of terms would have to 
be included to achieve an acceptable accuracy for the 
distribution function, which would then  
have to be numerically inverted (see their Section~7).

Our main new result and simulation method is based on 
the following theorem. The results proved rely
on the L\'evy characteristic function~\citep{Levy} for the L\'evy area 
$A(h)$; see Section~\ref{sec:proof}.
\begin{theorem}[Logistic Expansion]\label{th:logisticexpansion}
The L\'evy area $A(h)$ conditioned on the Wiener increments
$\Delta W^1$ and $\Delta W^2$ is equivalent in distribution
to the series of Logistic random variables
$A(h)\sim \lim_{N\to\infty}A_N(h)$, where
\begin{equation*}
A_N(h)\coloneqq
\frac{h}{2\pi}\biggl(X+\sum_{n=0}^N\frac{1}{2^n}\sum_{k=1}^{P_n}X_{n,k}\biggr),
\end{equation*}
where for $n=0,1,\ldots,N$: the 
$P_n\sim\textsf{Poisson}\bigl(\tfrac12 a^22^{n}\bigr)$ 
are independent Poisson random variables, 
for $k=1,2,\ldots,P_n$, $X=\log\bigl(U/(1-U)$ and
$X_{n,k}=\log\bigl(U_{n,k}/(1-U_{n,k})$ with  
$U,U_{n,k}\sim\textsf{Unif}([0,1])$ 
independent identically distributed 
uniform random variables
(i.e.\/ $X,X_{n,k}\sim\textsf{Logistic}(1)$ are
independent identically distributed Logistic random variables).
The mean-square error of the Logistic expansion approximation 
$A_N(h)$ is \emph{exactly} $a^2h^2/(3\cdot2^{N+3})$.
The Logistic approximation to $A(h)$ including simulating
the tail sum is 
\begin{equation*}
A_N(h)+(ah/\sqrt{3\cdot2^{N+3}})\,Z,
\end{equation*}
where $Z\sim\textsf{N}(0,1)$. The mean-square error in this 
approximation is bounded by $h^2/(15\cdot2^{2N+1})$.
\end{theorem}
In the Logistic approximation in the theorem, at each order $n$, 
we must on average sample and add $\mathbb E\,(P_n)=\frac12a^22^n$ 
Logistic random variables. 
The \emph{strong} mean-square error results imply that to achieve
accuracy $\eps$ we require $N$ such that $2^N$ is of order $h^2/\epsilon$.
Hence the complexity of the Logistic approximation $A_N(h)$ is $h^2/\eps$; 
with tail simulation it is $h/\eps^{1/2}$. However, if
we can simulate the sum of say $P\gg1$ independent Logistic random variables
efficiently, then our representation can be used as a basis for an effective
L\'evy area sampling method; see Section~\ref{sec:alg} for details. 
Indeed for $P=10^3,10^4,10^5,10^6$ we 
approximate the inverse distribution function for the sum
of $P$ independent Logistic random variables by Chebychev
polynomials, in the central, middle and tail regions of
the inverse distribution. With this replacement,
we still achieve a \emph{strong approximation}.
Our tail region stops at distance $10^{-12}$ from the endpoints. 
The error of the Chebychev approximations in the three regions 
is controlled by the degrees of the polynomials we prescribe. 
We choose to require uniform errors of order $10^{-12}$,
which is far smaller than the Monte Carlo error we could achieve, 
and which we regard as \emph{almost-exact}.  
Note to approximate the L\'evy area, we truncate
the Logistic series representation to include terms
with $n\leqslant N$. 
The mean-square truncation error implies that, as above,
we require $N$ such that $2^{N}$ is of order $h^2/\epsilon$.
However for our truncation we must on average sample order 
$N^2$ uniform random variables---this is the 
sum over $n\leqslant N$ of the sum of the digits in $\mathbb E\,(P_n)$.
Hence the complexity is the square of the logarithm of
$h^2/\epsilon$ without tail and square logarithm of $h/\epsilon^{1/2}$ with tail;
see the electronic supplement for more details.
To summarize, the advantages of our direct inversion method based
on the Logistic expansion are: (1) its square logarithmic complexity and 
(2) the main ingredient is direct inversion, which importantly, 
can be used in combination with quasi-Monte Carlo simulation.

\section{Direct inversion algorithm}
\label{sec:alg}
We apply the ideas underlying the Beasley--Springer--Moro method
for standard Normal random variables. 
More details, including comprehensive details of anaylsis and 
tables of polynomial coefficients, can be found in the electronic supplement.
We consider the fixed values $P=10^3,10^4,10^5,10^6$.
Let $\Phi$ denote the distribution function for the sum
of $P$ independent Logistic random variables. 
The inverse distribution function $\Phi^{-1}=\Phi^{-1}(u)$ 
is antisymmetric about $1/2$. We thus focus on the subinterval $[1/2,1)$
of its support. Indeed we split this interval into the three
regions: the \emph{central} $[1/2,u_1]$; \emph{middle} $(u_1,u_2]$ and 
\emph{tail} $(u_2,1-10^{-12}]$ \emph{regions}. 
We neglect the regions at distance $10^{-12}$ from the endpoints. 
The values $u_1=u_1(P)$ and $u_2=u_2(P)$ roughly separate
the characteristic behaviour of $\Phi^{-1}$. 
In the \emph{central region} we approximate $\Phi^{-1}\approx U\cdot C_n(z)$
where $C_n$ is a degree $n$ Chebychev polynomial approximation,
where $U\coloneqq(2P\pi^3/3)^{1/2}(u-1/2)$ and $z=k_1U^2+k_2$.
In the \emph{middle} and \emph{tail regions} 
we approximate $\Phi^{-1}\approx C_n(z)$ where
$z=k_1U+k_2$ and
$U\coloneqq\pi\bigl(-\frac23 P\log(2\sqrt{\pi}(1-u))\bigr)^{1/2}$.
Briefly the rationale underlying this ansatz for $U$ is as follows.
By the Central Limit Theorem, asymptotically in distribution, we have
$\Phi(x)\sim\Phi_{\mathrm{N}}\bigl(x/\pi\,(P/3)^{1/2}\bigr)$, where
$\Phi_{\mathrm{N}}$ is the standard Normal distribution
function. Following \citet{Moro}, using the asymptotic tail approximation
for the standard Normal, we find 
$\Phi(x)\sim 1-(1/2\pi^{1/2})\exp(-x^2/(2\pi^2P/3))$. 
Inverting this relation gives the ansatz for $U$.
Note, in all three regions, $k_1$ and $k_2$ are chosen to ensure
$z=-1$ at the left endpoint and $z=+1$ at the right endpoint.

\begin{figure}
  \begin{center}
  \includegraphics[width=16cm,height=4cm]{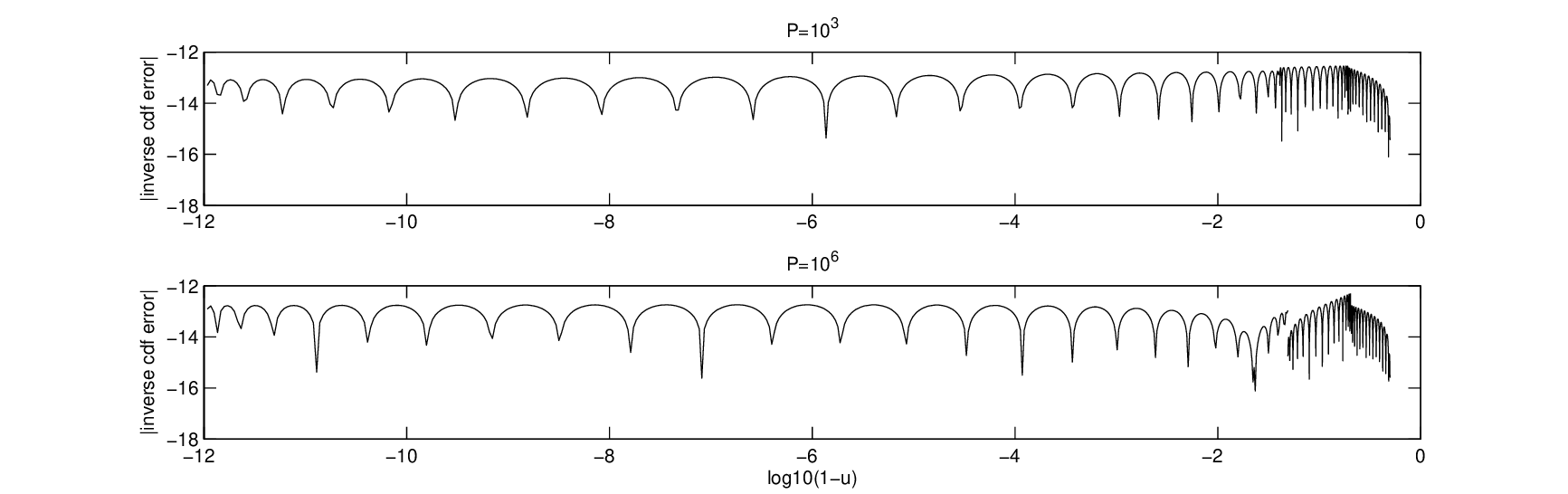}
  \end{center}
  \caption{The panels show the error in the Chebychev polynomial approximations
to $\Phi^{-1}=\Phi^{-1}(u)$ for $P=10^3$ (top) and $P=10^6$ (bottom) 
across all three regions with $u\in[1/2,1-10^{-12}]$. Note we use $1-u$ on the 
abscissa and a $\log$-$\log_{10}$ plot to highlight the tail region 
(on the left).}
\label{fig:cheberror}
\end{figure}

The coefficients for the Chebychev polynomial approximations
are computed in the standard fashion, see Section~5.8 in \citet{PTVF}.
However two additional aspects are crucial to their accurate and
efficient evaluation. First, based on the inverse 
Fourier transform of the characteristic function for the sum
of $P$ Logistic random variables, we derived a large $P$ asymptotic
approximation for the distribution function $\Phi=\Phi(x)$. 
Following \citet[pp.~272--4]{BO} we developed the
expansion in reciprocal powers of $P$ to all orders.
Indeed we computed up to $140$ terms to obtain the 
requisite accuracy for $\Phi$ in the tail regions. 
Second, we computed the expansion and performed the required
rootfinding for $\Phi^{-1}$ values using Maple with $25$ 
(and on occasion $45$) digit accuracy. We imported these
accurate and reliable Chebychev coefficients to Matlab. 
All subsequent computations are done using
double precision in Matlab---the Chebychev approximations were
evaluated using Clenshaw's recurrence formula \citep[p.~193]{PTVF}.
Figure~\ref{fig:cheberror} shows the errors in 
the Chebychev polynomial approximations for $\Phi^{-1}$
across all three regions when $P$ is $10^3$ and $10^6$.  
The Chebychev polynomials have degrees 
from $13$ to $15$ in the centre (roughly $[0.5,0.8]$) 
and middle regions (roughly $[0.8,0.95]$),
and degrees $24$ to $26$ in the tail region.

The direct inversion algorithm works as follows. We use the Logistic
series representation in Theorem~\ref{th:logisticexpansion} which
we truncate at some large integer $N$. First, consider the 
Poisson samples $P_n$ required for each $n=0,\ldots,N$. 
These are obtained by direct inversion
when the mean $\frac12a^22^n$ is $100$ or less---rather
than faster acceptance-rejection methods.
For means $100$ or greater, we use the PTRS transformed rejection 
(or almost-exact inversion) method from \citet{H}. 
To achieve an almost-exact 
quasi-Monte Carlo implementation, tables of the Poisson distribution
function could be constructed for different representative means
such as $10^i$ for $3\leqslant i\leqslant12$. Since the standard
deviation of a Poisson random variable 
is the square-root of the mean, such tables which ignore
tails of order $10^{-12}$ will not be restrictively large. 
Samples can be drawn for these representative 
means by a fast look-up algorithm. A Poisson sample for a given mean 
could be generated by adding the requisite numbers of Poisson samples 
from the representative means. Second, the sums of $P_n$ Logistic
random variables are handled thus. Whenever $P_n<10^3$ we 
sum $P_n$ Logistic random variables. However if $P_n\geqslant10^3$
we decompose $P_n=p+p_3\cdot10^3+p_4\cdot10^4+p_5\cdot10^5+p_6\cdot10^6$.
Here $p=P_n\,\text{mod}\,10^3$ and the $p_k$ are the multiples of $10^k$
present in the sample $P_n$; note $p_6$ need not be a single digit.
For each $k=3,4,5,6$ we sample $p_k$ random variables from the distribution
function for the sum of $10^k$ Logistic random variables
using the corresponding Chebychev polynomial approximations described above.

\section{Simulations}\label{sec:sim}
\begin{figure}
  \begin{center}
  \includegraphics[width=16cm,height=6cm]{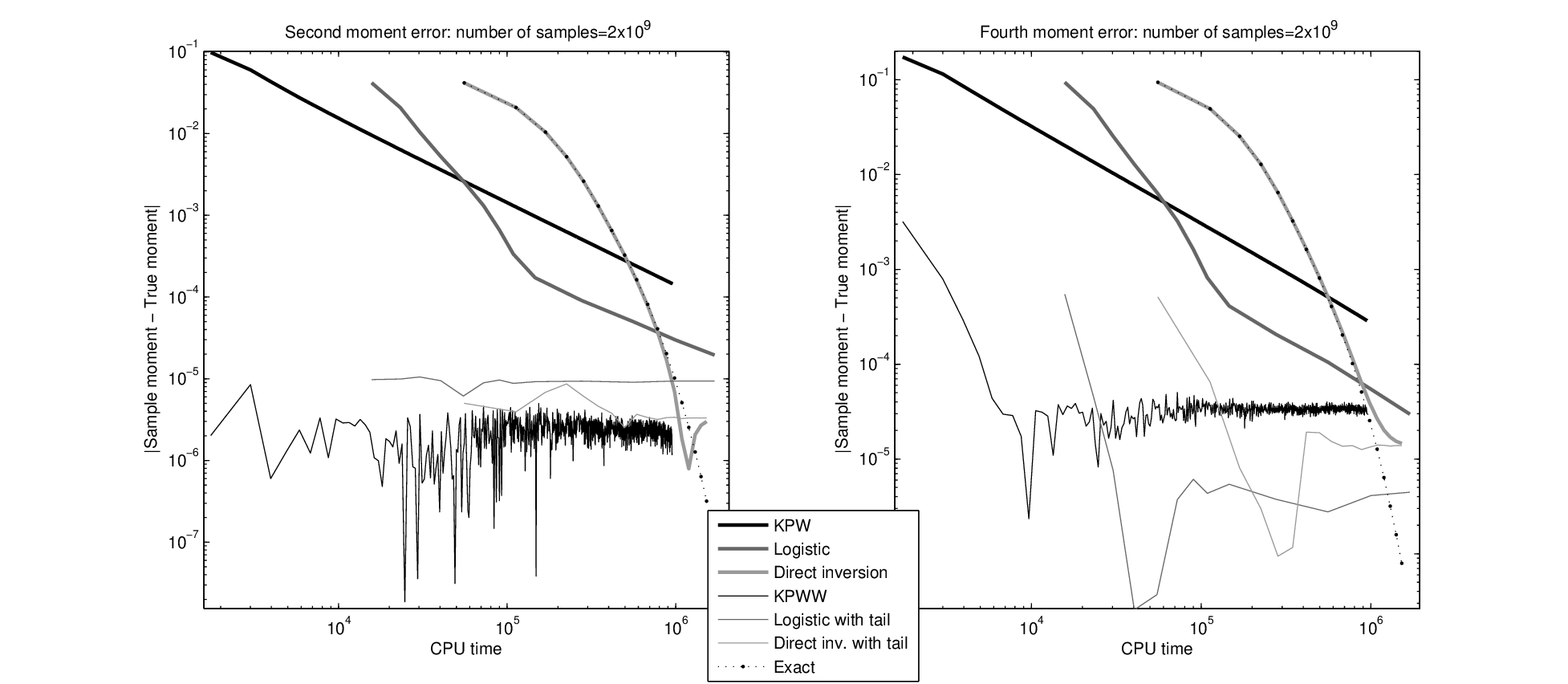}
  \end{center}
  \caption{We show the absolute error in the second (left) and fourth (right) moments
versus the CPU time required to compute the simulation for three different
methods (with and without tail simulation). The methods are: 
Kloeden--Platen--Wright (KPW); Logistic method based on adding the requisite
Logistic random variables; direct inversion based on the Logistic
expansion; Kloeden--Platen--Wright--Wiktorsson (KPWW); 
Logistic expansion with tail simulation and direct inversion 
with tail simulation. We fixed $h=1$, hence the exact second and fourth moments 
are $1/4$ and $5/16$, respectively. The corresponding exact error in these
moments for the direct inversion method is also indicated.}
\label{fig:difference}
\end{figure}

We simulated the L\'evy area $A(h)$ using three methods with and 
without tail simulation for $h=1$. These are: 
Kloeden--Platen--Wright from the introduction (no tail); 
Logistic method using Theorem~\ref{th:logisticexpansion} where
the requisite Poisson number of Logistic random variables are added
at each order $n=0,\ldots,N$ and direct inversion based on the Logistic expansion 
as described in Section~\ref{sec:alg}. We also implemented these three methods with
tail simulation as shown in the introduction and
Theorem~\ref{th:logisticexpansion}, respectively.
In Figure~\ref{fig:difference} the panels show the absolute error in the second moment
(left) and fourth moment (right)
versus the CPU time required to compute the simulation.
Since $\mathbb E\,a^2=\mathbb E\,\bigl((\Delta W^1)^2+(\Delta W^2)^2\bigr)/h=2$ and $h=1$,
the true variance $(1+\mathbb E\,a^2)\,h^2/12$ is $1/4$. 
The exact fourth moment is $5/16$. 
We performed $L=2\times10^9$ simulations in all cases. The Monte Carlo error
is of order $L^{-1/2}$. This can be observed in Figure~\ref{fig:difference}, 
where the error curves become ``horizontal and noisy''. For all
methods we truncated after $N$ terms, increasing $N$ in integers for 
the Kloeden--Platen--Wright method, and in powers of $2$ for the
Logistic and direct inversion methods. We included
$N=2^{10},2^{13}$ and $2^{18}$ terms for the Kloeden--Platen--Wright, 
Logistic and direct inversion methods, respectively.
The errors of the sample moments in the basic approximations 
in Figure~\ref{fig:difference} are in good correspondence
with their theoretical values. We observe this for example, 
for the direct inversion method for which we also plot the exact 
second and fourth moment errors versus the corresponding
simulation CPU times (the fourth moment error calculation can
be found in the electronic supplement). In particular we observe
the behaviour of the exact error in the large $N$ regime, where
the simulations are dominated by the Monte Carlo error. 
We also observe the complexity of $\varepsilon^{-1}$ for the Kloeden--Platen--Wright
and Logistic expansion methods (without tail simulation).
In the Logistic method case this complexity is
observed in the high accuracy/large effort asymptotic limit.
Although the direct inversion method requires more effort at
low accuracies, the square logarithmic complexity is 
observed in the high accuracy/large effort asymptotic limit.
In each method the tail approximations are designed
so that the tail is approximated by a matched Normal random variable
(the first two moments are matched). Hence they simulate the 
second moment exactly and we only observe Monte Carlo noise
on the left in Figure~\ref{fig:difference}.
The simulations thus confirm the analysis and expected properties.

\section{Proof of the Logistic Expansion Theorem}\label{sec:proof}
The \emph{characteristic function} $\hat\phi=\hat\phi(\xi)$ 
corresponding to the probability density
function $\phi$ for the L\'evy area $A(h)$, 
given $\Delta W^1$ and $\Delta W^2$, is
$\hat\phi(\xi)=\bigl(\tfrac12 h\xi/\sinh(\tfrac12h\xi)\bigr)
\exp\bigl(-\tfrac12a^2\bigl(\tfrac12h\xi\coth(\tfrac12h\xi)-1\bigr)\bigr)$.
We observe the identity $\mathrm{coth}\,z\equiv\mathrm{coth}\,z/2-1/\mathrm{sinh}\,z$,
iterated $N$ times, generates the identity 
$\mathrm{coth}\,z\equiv\mathrm{coth}(z/2^{N+1})-\sum_{n=0}^N1/\mathrm{sinh}(z/2^n)$.
Substituting this identity into the characteristic function we see that
\begin{equation*}
\exp\bigl(-\tfrac12a^2(z\,\mathrm{coth}\,z-1)\bigr)
=\prod_{n=0}^N\exp\biggl(\tfrac12a^22^n\cdot\frac{z/2^n}{\mathrm{sinh}\,z/2^n}\biggl)
\exp\bigl(-\tfrac12a^22^n\bigr)\cdot\exp\bigl(\mathcal E_N(z,a)\bigr),
\end{equation*}
where 
$\mathcal E_N(z,a)\coloneqq 
-\tfrac12a^2\bigl(z\,\mathrm{coth}\,(z/2^{N+1})-2^{N+1}\bigr)\to0$ 
as $N\to\infty$, $\forall\,z\in\mathbb R$.
Thus $\forall\,z=\xi h/2\in\mathbb R$, we have
\begin{equation*}
\hat\phi(\xi)=\frac{\xi h/2}{\mathrm{sinh}\,\xi h/2}
\cdot\prod_{n=0}^\infty \mathbb E\,
\biggl(\frac{\xi h/2^{n+1}}{\mathrm{sinh}\,\xi h/2^{n+1}}\biggr)^{P_n},
\end{equation*}
where $P_n\sim\textsf{Poisson}\bigl(\tfrac12 a^22^{n}\bigr)$.
Note that the expression
$(\xi h/2^{n+1})/\mathrm{sinh}(\xi h/2^{n+1})$ is the 
characteristic function of a
$\textsf{Logistic}(h/2^{n+1}\pi)\sim (h/2^{n+1}\pi)\cdot\textsf{Logistic}(1)$
random variable. The Logistic expansion follows. 

We now focus on the error statements. First consider the tail sum itself. 
Directly computing
\begin{equation*}
\mathbb E\,\bigl|A(h)-A_N(h)\bigr|^2=
\mathbb E\,\biggl(\frac{h}{2\pi}\sum_{n=N+1}^\infty 
\frac{1}{2^{n}} \sum_{k=1}^{P_n}
X_{n,k}\biggr)^2
=\Bigl(\frac{h}{2\pi}\Bigr)^2\sum_{n=N+1}^\infty 
\frac{1}{2^{2n}} \sum_{l=0}^\infty\mathbb P\,\{P_n=l\}
\cdot\mathbb E\,\biggl(\sum_{k=1}^{l} X_{n,k}\biggr)^2.
\end{equation*}
Now we use that the expectation on the far right equals $\ell\pi^2/3$ and
that $\sum_{l\geqslant0}\mathbb P\,\{P_n=l\}\cdot l=\mathbb E\,\{P_n\}=a^22^{n-1}$.
Then noting that $\sum_{n\geqslant N+1}2^{-2n}=2^{-(N+1)}$ gives the exact
mean-square error result. 

Second, we derive the error bound for the approximation including the 
Normal tail sum approximation. The tail sum is an infinitely divisible class G
random variable, i.e.\/ its characteristic function has the form $\hat{\phi}(\xi) =
\exp(-\Psi(\xi^2))$, where $\Psi(0)=0$, and $(-1)^{n-1} \Psi^{(n)}(\xi)\ge
0$ for all $n$, see \citet{RW} (pg.~163). Using that
$A(h)-A_N(h)=\frac{h}{2\pi}\sum_{n\geqslant N+1}2^{-n}\sum_{k=1}^{P_n}X_{n,k}
=\frac{\tilde{h}}{2\pi}\sum_{n\geqslant0}2^{-n} \sum_{k=1}^{Q_{n}}X_{n,k}$,
where $\tilde{h}=h/2^{N+1}$, and where $Q_n$ has a Poisson
distribution with parameter $\frac{1}{2} \tilde{a}^2 2^n$ and where $\tilde{a}^2=a^2
2^{N+1}$, it follows that
the tail sum $A(h)-A_N(h)$ has the characteristic function \citep{Levy}
$\hat{\phi}_N(\xi)=\exp\bigl(-\tfrac12\tilde{a}^2\bigl(\tfrac12\tilde{h}
\xi\coth(\tfrac12\tilde{h}\xi)-1\bigr)\bigr)
=\exp\bigl(-\tfrac12\tilde{a}^2
\sum_{n=1}^\infty\xi^2/\bigl((n\tilde{h}/2\pi)^2+\xi^2\bigr)\bigr)$.
Hence $A(h)-A_N(h)$ has class G distribution---see also Prop.~5 in 
\citet{RW}. We can now proceed as in the proof of 
Theorem~7 in Ryd\'en and Wiktorsson. The tail sum  can be represented as a
product of a standard Normal
random variable $Z$ and the square root of an independent 
positive, infinitely divisible variable
$Y_N$, i.e.\/ $A(h)-A_N(h)=Z\sqrt{Y_N}$.
If $\sigma_N^2$ denotes the variance of $A(h)-A_N(h)$,
then the mean-square error when including the Normal tail
approximation is given by
\begin{equation*}
\mathbb E\,\bigl|A(h)-A_N(h)-\sigma_N Z\bigr|^2
=\mathbb E\,\{Z^2\}\cdot \mathbb E\,\bigl|\sqrt{Y_N}-\sigma_N\bigr|^2 
=\mathbb E\,(Y_N-\sigma_N^2)^2/(\sqrt{Y_N}+\sigma_N)^2.
\end{equation*}
This is bounded above by $\mathbb E\,\bigl(Y_N-\sigma_N^2\bigr)^2/\sigma_N^2$.
Let $g$ denote the Laplace transform of $Y_N$. Then
$g(z)=\hat{\phi}_N(\sqrt{2z})$, and
the variance of $Y_N$ is given by $(\log g)''(0)$.
If $\ell(z)\coloneqq1-\sqrt{2z}\coth(\sqrt{2z})$
then we see that 
$\log g(z)=\tilde{a}^2 \ell(z\tilde{h}^2/4)$
and $(\log g)''(0)=\tilde{a}^2(\tilde{h}/2)^4 \ell{''}(0)$.
Hence the mean-square error when including tail
approximation is bounded by 
$\sigma_N^{-2}\tilde{a}^2(\tilde{h}/2)^4\ell''(0)
=h^2/(15\cdot2^{2N+1})$, completing the proof.

\section{Conclusion}\label{sec:conclu}
We conclude with some brief observations. For low to medium accuracies
the Kloeden--Platen--Wright--Wiktorsson and Logistic expansion methods
perform extremely well. When high accuracies are
required and/or a quasi-Monte Carlo implementation is intended, then
the almost-exact direct inversion method is the method of choice
as can be observed in Figure~\ref{fig:difference}. In the electronic
supplement we provide more detailed complexity calculations
in order to establish a cross-over accuracy criterion. 
In other words we specify the accuracy for which one would 
preferentially use the direct inversion over the 
Kloeden--Platen--Wright--Wiktorsson method. This depends of course
on an efficient implmentation of the algorithm and the system
on which it is implemented. Though we endeavoured to establish
the three regions and Chebychev polynomial approximations for each $P$
so that the direct inversion algorithm is efficient, they could be 
further optimized. Note for example, for our simulation implementation in
Figure~\ref{fig:difference} for medium accuracies we would preferentially
use the basic Logistic method.
The direct inversion techniques above could also be applied to the 
Ryd\'en--Wiktorsson series which involves large sums of independent 
Laplace random variables. 
Lastly, we remark on the case of a $d$-dimensional Wiener process
with $d\geqslant3$. When $d=3$, the characteristic variable is
$Z\in\mathfrak{so}(3)$ and $Z^3=-\|\zb\|^2\,Z$,
where $\|\zb\|$ is the Euclidean norm of the vector $\zb$ of the
three upper triangular components of $Z$. This property simplifies
power series functions of $Z$ with scalar coefficients and underlies 
Rodrigues formulae for $Z\in\mathfrak{so}(3)$. Such formulae 
generalize to higher dimensions; see \citet{GX}. In particular,
the joint characteristic function given by \citet{W} for the 
L\'evy areas $A_{12}$, $A_{13}$ and $A_{23}$ conditioned on 
$\Delta\boldsymbol W\coloneqq(\Delta W^1,\Delta W^2,\Delta W^3)^{\text{\tiny{T}}}$
reduces to (scaling $\zb$ by $h/2$)
\begin{equation*}
\frac{\|\zb\|}{\sinh\|\zb\|}
\cdot\exp\bigl(-\tfrac12\|\Delta\boldsymbol W\|^2\,(\|\zb\|\,\mathrm{coth}\|\zb\|-1)\bigr)
\cdot\exp\bigl(\tfrac12\|\Delta\boldsymbol W\|^2\,
\langle\Delta\hat{\boldsymbol W},\hat\zb\rangle^2\,
(\|\zb\|\,\mathrm{coth}\|\zb\|-1)\bigr).
\end{equation*}
Here $\langle\Delta\hat{\boldsymbol W},\hat\zb\rangle$ is the Euclidean inner product
of the corresponding unit vectors. This result can be found in \citet[p.~23]{MY}.
The first factor is the characteristic function of a generalized Logistic random
variable and the second (radial-type) factor can be analysed in the same fashion 
as we have presented here. Appropriately efficiently simulating the third (angular-type) 
factor is our next goal.

\section*{Acknowledgements}
The authors thank the referee for their valuable and particularly helpful suggestions,
in particular for suggesting we include the exact errors in Figure~\ref{fig:difference}
and compute a ``cross-over'' accuracy.

\bibliographystyle{elsarticle-harv}

\newpage

\section*{Supplementary material}
We present comprehensive details underlying some of the results 
in our manuscript. We include:  (1) a derivation of the exact
error in the fourth moment; (2) an accuracy criterion which indicates
when, in terms of complexity, it is more preferential to use the direct inversion 
over Kloeden--Platen--Wright--Wiktorsson method; and for a sum of Logistic 
random variables, we include (3) the derivation of the distribution function 
representation in the large sum asymptotic limit; 
(4) complete tables of the coefficients used for the Chebychev polynomial 
approximations for the inverse distribution function and
(5) the derivation of a finite representation for the density function.

\subsection*{(1) Fourth moment error}\label{sec:fourthmomenterror}
In the main text we plotted the exact error in the fourth moment of the 
Logistic expansion versus CPU time from simulations for the
direct inversion method---see Figure~2 (right panel). 
The analytical expression for the exact error is as follows.
\begin{lemma}
The error in the fourth moment of the approximation $A_N(h)$
to the L\'evy area $A(h)$, 
given $a^2\coloneqq\bigl((\Delta W^1)^2+(\Delta W^2)^2\bigr)/h$, is
\begin{equation*}
\mathbb E\,\bigl(A_N(h)^4\bigr)-\mathbb E\,\bigl(A(h)^4\bigr)
=\frac{a^2h^2}{2^{N+4}}\cdot\biggl(-\frac13(1+a^2)+\frac{a^2}{12\cdot 2^N}
-\frac{1}{30\cdot 4^N}\biggr).
\end{equation*}
\end{lemma}
\begin{proof}
First we note that the fourth moment of the L\'evy area is given by the 
fourth derivative of the characteristic function evaluated at the origin.
This can be straightforwardly be computed and is given by
\begin{equation*}
\mathbb E\,\bigl(A(h)^4\bigr)
=\frac{1}{2^4}\Bigl(\frac{7}{15}+\frac{14}{15}a^2+\frac{1}{3}a^4\Bigr).
\end{equation*}
Second we compute $\mathbb E\,\bigl(A_N(h)^4\bigr)$. For convenience 
we set 
\begin{equation*}
Y_0\coloneqq X\qquad\text{and}\qquad 
Y_{n+1}\coloneqq\frac{1}{2^n}\sum_{k=1}^{P_n}X_{n,k}
\end{equation*}
for  $n=0,1,\ldots,N$, so that $A_N(h)$ can be expressed in the form
\begin{equation*}
A_N(h)=\frac{h}{2\pi}\sum_{n=0}^{N+1}Y_n.
\end{equation*}
Here, as stated in the Logistic Expansion Theorem in the main text, 
the $X$ and $X_{n,k}$ 
are independent identically distributed Logistic random variables and
the $P_n$ are independent Poisson random variables with expectation
$\frac12a^22^n$. Hence explicitly we have
\begin{align*}
\mathbb E\,\bigl(A_N(h)^4\bigr)
=&\;\Bigl(\frac{h}{2\pi}\Bigr)^4
\cdot\mathbb E\,\biggl(\sum_{n=0}^{N+1}Y_n\biggr)^4\\
=&\;\Bigl(\frac{h}{2\pi}\Bigr)^4
\cdot\biggl(\sum_{n=0}^{N+1}\mathbb E\,Y_n^4+
6\sum_{n=0}^{N}\sum_{m=n+1}^{N+1}\mathbb E\,Y_n^2\cdot\mathbb E\,Y_m^2\biggr).
\end{align*}
Note that $\mathbb E\,Y_0^2=\pi^2/3$ and $\mathbb E\,Y_0^4=7\pi^4/15$. 
For $n=1,\ldots,N+1$, similar to the error statement calculations in the
proof of the Logistic Expansion Theorem in Section~4 of the main text, we have
\begin{align*}
\mathbb E\,Y_n^2&\;=\frac{1}{2^{2(n-1)}}\cdot
\mathbb E\,\biggl(\sum_{k=1}^{P_{n-1}}X_{n-1,k}\biggr)^2\\
&\;=\frac{1}{2^{2(n-1)}}\cdot\mathbb E\,\{P_{n-1}\}
\cdot\frac{\pi^2}{3}\\
&\;=\frac{1}{2^{2(n-1)}}\cdot\frac12a^22^{n-1}\cdot\frac{\pi^2}{3}.
\end{align*}
Further we can similarly compute for $n=1,\ldots,N+1$: 
\begin{align*}
\mathbb E\,Y_n^4&\;=\frac{1}{2^{4(n-1)}}\cdot
\mathbb E\,\biggl(\sum_{k=1}^{P_{n-1}}X_{n-1,k}\biggr)^4\\
&\;=\frac{1}{2^{4(n-1)}}\cdot
\sum_{\ell=0}^\infty \mathbb P\,\{P_{n-1}=\ell\}
\cdot\mathbb E\,\biggl(\sum_{k=1}^{\ell}X_{n-1,k}\biggr)^4\\
&\;=\frac{1}{2^{4(n-1)}}\cdot
\sum_{\ell=0}^\infty \mathbb P\,\{P_{n-1}=\ell\}
\cdot\Bigl(\frac{7\pi^4}{15}\ell+6\cdot\frac{\pi^4}{9}\cdot\frac12\ell(\ell-1)\Bigr)\\
&\;=\frac{1}{2^{4(n-1)}}\frac{\pi^4}{3}\cdot
\sum_{\ell=0}^\infty \mathbb P\,\{P_{n-1}=\ell\}
\cdot\Bigl(\ell^2+\frac25\ell\Bigr)\\
&\;=\frac{1}{2^{4(n-1)}}\frac{\pi^4}{3}\cdot
\mathbb E\,\Bigl\{P_{n-1}^2+\frac25P_{n-1}\Bigr\}\\
&\;=\frac{1}{2^{4(n-1)}}\frac{\pi^4}{3}\cdot
\Bigl(\Bigl(\frac12a^22^{n-1}+\frac14a^42^{2(n-1)}\Bigr)+\frac25\cdot\frac12a^22^{n-1}\Bigr)\\
&\;=\frac{1}{2^{4(n-1)}}\frac{\pi^4}{3}\cdot\Bigl(\frac14a^42^{2(n-1)})+\frac{7}{10}a^22^{n-1}\Bigr).
\end{align*}
Using these last two results in our expression for $\mathbb E\,\bigl(A_N(h)^4\bigr)$
above we find 
\begin{align*}
\mathbb E\,\bigl(A_N(h)^4\bigr)=&\;
\Bigl(\frac{h}{2}\Bigr)^4\cdot\biggl(\frac{7}{15}
+\sum_{n=1}^{N+1}\frac13\Bigl(\frac14a^42^{-2(n-1)})+\frac{7}{10}a^22^{-3(n-1)}\Bigr)\\
&\;\qquad\qquad\quad\;+6\sum_{m=1}^{N+1}\frac13\cdot\frac13\cdot\frac12a^22^{-(m-1)}
+6\sum_{n=1}^{N}\sum_{m=2}^{N+1}\frac19\cdot\frac14a^42^{-(n+m-2)}\biggr)\\
=&\;\Bigl(\frac{h}{2}\Bigr)^4\cdot\biggl(\frac{7}{15}
+a^2\Bigl(\frac{14}{15}-\frac232^{-N-1}+\frac{4}{15}8^{-N-1}\Bigr)
+a^4\Bigl(\frac13-\frac132^{-N}+\frac134^{-N-1}\Bigr)\biggr),
\end{align*}
where we used the formula for partial geometric sums.
Taking the difference of this result with the expression above for 
$\mathbb E\,\bigl(A(h)^4\bigr)$ we establish the result.
\end{proof}
\begin{remark}
Note that in Figure~2 we take $h=1$ and use that
$\mathbb E\,a^2=2$ and  $\mathbb E\,a^4=8$.
\end{remark}
\begin{remark}
Using analogous arguments to those in the proof of the Logistic
Expansion Theorem in Section~4 of the main text and some
of the results above, it is straightforward to derive an exact
analytical expression for $\mathbb E\,|A(h)-A_N(h)|^4$. We can also
derive an upper bound for $\mathbb E\,|A(h)-A_N(h)-\sigma_NZ|^4$,
where $\sigma_NZ$ represents the approximation to the tail sum 
(as in the main text).
\end{remark}

\subsection*{(2) Accuracy criterion}\label{sec:accuracycriterion}
At what accuracy is the complexity of the direct inversion
method based on the Logistic expansion better than the 
Kloeden--Platen--Wright--Wiktorsson method? As mentioned 
before this depends on how the algorithms are implemented
and what system they are implemented on.
We use the number of uniform random variables required 
to generate a L\'evy area sample on each timestep
to achieve a mean-square accuracy of $\eps$, as the measure
of complexity. Consider the Kloeden--Platen--Wright method.
At leading order the mean-square error, when we truncate
the series to include the terms $n\leqslant N$, is equal to 
$C_{\mathrm{KPW}}\cdot h^2/N$,
where $C_{\mathrm{KPW}}$ is a known constant.
Equating the mean-square error
with $\eps$, we see that to achieve a mean-square accuracy of $\eps$ 
we must choose $N=C_{\mathrm{KPW}}\cdot h^2/\eps$.
The computational effort required to generate the truncated series
is at leading order given by $A_{\mathrm{KPW}}\cdot N$,
where $A_{\mathrm{KPW}}$ is a positive constant. Indeed at each order, 
for each term we need to generate $4$ uniform random variables 
(for each of the $4$ standard normal random variables required) 
and we need to also account for the flops and function evaluations 
required to compute the term. We absorb these factors multiplicatively into the 
constant $A_{\mathrm{KPW}}$. 
Hence the Kloeden--Platen--Wright method \emph{complexity} is 
\begin{equation*}
B_{\mathrm{KPW}}\cdot h^2/\eps,
\end{equation*}
where we set $B_{\mathrm{KPW}}\coloneqq A_{\mathrm{KPW}}\cdot C_{\mathrm{KPW}}$.

For the direct inversion method based on the Logistic expansion
the mean-square error at leading order when we truncate to 
include the terms $n\leqslant N$, is equal to
$C_{\mathrm{DI}}\cdot h^2/2^N$,
where $C_{\mathrm{DI}}=a^2/24$. Again, equating the mean-square error
with $\eps$, to achieve a mean-square accuracy of $\eps$ 
we must choose $N=\log(a^2h^2/24\,\eps)/\log 2$.
We postulate that the computational effort required to generate this approximation
for the direct inversion method has the general quadratic form
\begin{equation*}
A+BN+CN^2,
\end{equation*}
for some constants $A$, $B$ and $C$ dependent on the implementation.
Hence the \emph{complexity} of the direct inversion method has this
form with $N=\log(a^2h^2/24\,\eps)$, absorbing $1/\log2$ factors into
the constants. The rationale for this quadratic form is as follows.
For each $n\leqslant N$ we generate a Poisson random variable 
$P_n$ whose expectation is proportional to $a^22^n$. 
Using direction inversion for the 
distribution functions of large sums of independent identically
distributed Poisson random variables, the number of uniform
random variables we need to generate is proportional to the
sum of digits in $P_n$ which is proportional
to $n$. Hence for all the terms $n\leqslant N$, the total
number of uniform random variables we need to generate 
is proportional to $N(N-1)$. Further, in Figure~2 in the main
text we observe for small values of $N$, there is
a non-negligible additive shift in the direct inversion complexity
compared to the basic Logistic method. This is largely
due to the algorithm for splitting the Poisson random
variable into its digits (which requires further optimization).
Hence we include $N^1$ and $N^0$ dependent terms in our complexity
form.

The direct inversion complexity is less
than the Kloeden--Platen--Wright complexity when
\begin{equation*}
A+B\cdot\bigl(\log(a^2h^2/24\,\eps)\bigr)+C\cdot\bigl(\log(a^2h^2/24\,\eps)\bigr)^2
<B_{\mathrm{KPW}}\cdot h^2/\eps.
\end{equation*}
All the quantities involved are known with the exception
of the constants $B_{\mathrm{KPW}}$, $A$, $B$ and $C$. These
encode the flops and function evaluations required
to compute the direct inversion and Kloeden--Platen--Wright
samples. They depend on the efficiency of the algorithms
implemented, programming language used and the hardware utilized.
The Kloeden--Platen--Wright method is very simple and these
factors are unlikely to influence $B_{\mathrm{KPW}}$ 
(indeed they do not; see below). However
as we mentioned in the Conclusion of the main text and above, 
the implementation of the direct inversion method we advocate
can be further optimized (though we endeavoured to implement 
it as efficiently as we could using compiled Matlab code).
The constants importantly affect the crossover accuracy below 
which the inequality above holds. 

\begin{figure}
  \begin{center}
  \includegraphics[width=10cm,height=8cm]{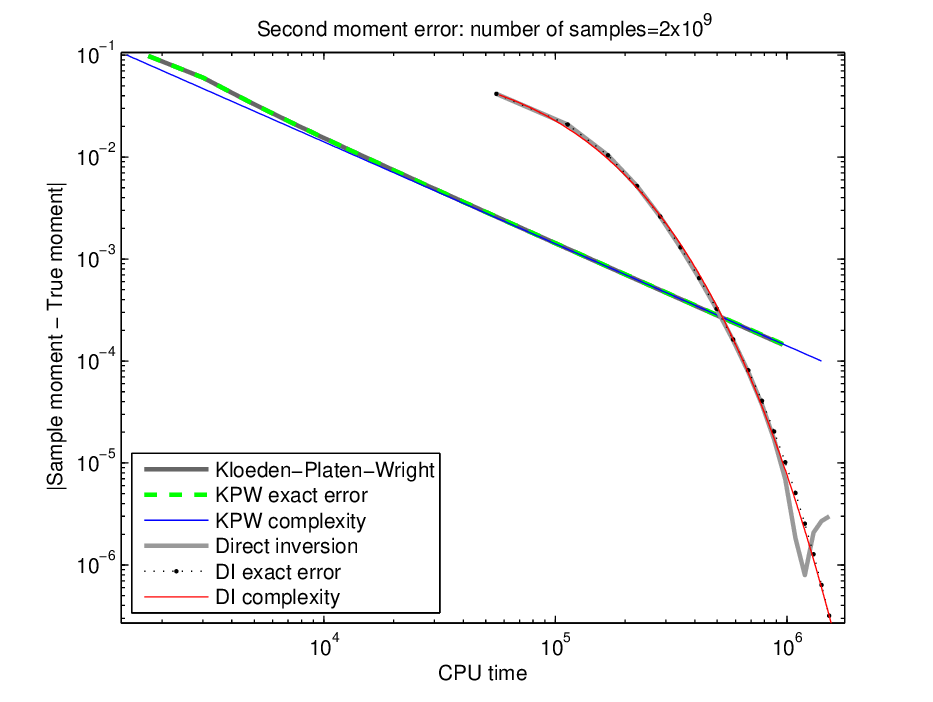}
  \end{center}
  \caption{We show the absolute error in the second moment
versus the CPU time required to compute the simulation for the 
Kloeden--Platen--Wright method (KPW) and direct inversion method 
based on the Logistic expansion (DI), as well as the exact error 
for the DI method---these are the same as in Figure~2 in the main
text. However in addition we include here the exact error for the
KPW method and the expressions for the complexity at
leading order (abscissa) plotted versus accuracy (ordinate) for the 
KPW method, $B_{\mathrm{KPW}}/\eps$, and DI method, $A+BN+CN^2$ 
with $N=\log(1/12\,\eps)$, derived in the text---the constants 
$B_{\mathrm{KPW}}$, $A$, $B$ and $C$ are fitted as described therein.
In our implementation the crossover accuracy is $10^{-3.5728}$.}
\label{fig:crossover}
\end{figure}

Let us estimate the crossover accuracy
in the context of our implementation generating Figure~2
in the main text---note the L\'evy area and the approximations
considered here have zero mean and are sums of independent random
variables, so the mean-square error and difference
of their second moments are the same. 
Recall that we set $h=1$ and that $\mathbb E\,a^2=2$.
The log-log plot in Figure~2 reveals the $B_{\mathrm{KPW}}/\eps$
behaviour we expect for the Kloeden--Platen--Wright method.
The mean-square error, from the analytical formula in the 
Introduction in the main text, for the case $N=2^9$ is 
$2.9655\times 10^{-4}$. The corresponding CPU time from 
our simulations is $4.7533\times10^5$.
Hence we have $B_{\mathrm{KPW}}=1.4096\times10^2$.
In Figure~\ref{fig:crossover} in this supplement, 
we plot accuracy $\eps$ versus $B_{\mathrm{KPW}}/\eps$
for this value of $B_{\mathrm{KPW}}$ and note
a markedly close fit to the simulation error
versus CPU time plot. The exact mean-square error 
plotted versus the simulation CPU times very
closes fits the simulation error
versus CPU time plot. There is thus very good
agreement between the simulations and our 
leading order complexity expression for 
the Kloeden--Platen--Wright method.
For the direct inversion method we match the
quadratic logarithmic form $A+BN+CN^2$ with
$N=\log(1/12\,\eps)$ to the three simulation
data points corresponding to $N$ equal to $1$, $9$
and $18$, where respectively the exact accuracy 
$\eps$ is $2^{-1}/12$, $2^{-9}/12$ and $2^{-18}/12$
and the corresponding CPU times are 
$5.5585\times10^4$, $5.8620\times10^5$ and $1.5223\times10^6$. 
With these values we deduce $A=9.2103\times10^3$,
$B=1.4669\times10^5$ and $C=2.4465\times10^4$.
We see in Figure~\ref{fig:crossover} in this supplement,
this fitted quadratic logarithmic form plotted
with complexity on the abscissa and accuracy on the 
ordinate, fits the error versus CPU time plot very
closely. In particular, we have established direct
practical confirmation of our main complexity claims
for the direct inversion method. Further we can 
compute the smallest root of $B_{\mathrm{KPW}}/\eps-(A+BN+CN^2)$ 
with $N=\log(1/12\,\eps)$, which is $10^{-3.5728}$
(which can be observed in Figure~\ref{fig:crossover}). 
For accuracies below this value, for our implementation,
the direct inversion method is more efficient than
the Kloeden--Platen--Wright method. With further 
optimization of the direct inversion implementation,
we expect this value to increase.

The computation efforts for the Kloeden--Platen--Wright--Wiktorsson
and direct inversion with tail methods are at leading order the 
same as for the Kloeden--Platen--Wright and direct inversion methods
as outlined above (only one further random variable needs to be sampled).
The accuracies are better: at leading order proportional to 
$h^2/N^2$ and $h^2/2^{2N}$, respectively. However currently, 
we only have upper bounds for the coefficients of these 
leading order mean-square errors, so we cannot make 
a meaningful direct comparison of complexities as above.  
In practice, at this stage in development, we suggest
that the computational effort corresponding to the crossover 
value for the Kloeden--Platen--Wright and direct inversion methods
gives a rough guide of the crossover computational effort value 
for the Kloeden--Platen--Wright--Wiktorsson
and direct inversion with tail methods. 

\subsection*{(3) Large sums of Logistic random variables}\label{sec:largesum}
We derive the asymptotic series representation for the distribution
function for the sum of $P$ Logistic random variables in the limit
$P\to\infty$. This series is crucial to computing the coefficients 
in our Chebychev polynomial approximations to the inverse distribution
function across its entire domain, efficiently and accurately.
Our asymptotic series expansion is developed using 
the inverse Fourier transform of the characteristic function for the sum 
of $P$ Logistic random variables. We apply the standard Laplace
method techniques outlined in \citet[p.~272--3]{BO} to
derive all the higher order terms in reciprocal powers of $P$.
Using the characteristic function for the sum of $P$
Logistic random variables, the probability density function as the 
inverse Fourier transform can be expressed in the form 
\begin{equation*}
\phi(x)\coloneqq\frac{1}{2\pi^2}\int_{\mathbb{R}}
\biggl(\frac{z}{\mathrm{sinh}\,z}\biggr)^P\cos\Bigl(\frac{xz}{\pi}\Bigr)\,\mathrm{d}z.
\end{equation*}
Using this form we prove the following representation
for the corresponding distribution function $\Phi=\Phi(x)$.
\begin{theorem}
The distribution function $\Phi=\Phi(x)$ has the asymptotic series
expansion as $P\to\infty$, 
\begin{equation*}
\Phi(x)\sim\frac{1}{2}+\biggl(\frac{3}{2\pi^3P}\biggr)^{1/2}
\biggl(\sum_{k\geqslant0}(-1)^k\frac{a_k}{(2k+1)!}\frac{x^{2k+1}}{P^k}
+\sum_{k\geqslant2}\sum_{\ell=2}^{k}\sum_{j=1}^{\llcorner\ell/2\lrcorner}
\frac{(-1)^{k-\ell}}{\pi^{2(k-\ell)}}\frac{a_k\cdot c_{j,\ell}}{\bigl(2(k-\ell)+1\bigr)!j!}
\frac{x^{2(k-\ell)+1}}{P^{k-j}}\biggr).
\end{equation*}
Here the constants $a_k$ and $c_{\ell,j}$ are given by 
\begin{equation*}
a_k\coloneqq3^k\cdot(2k-1)(2k-3)\cdots(1)
\qquad\text{and}\qquad 
c_{j,\ell}\coloneqq\sum_{i_1+\cdots+i_j=\ell-2j}\hat\varphi_{i_1}\cdots\hat\varphi_{i_j},
\end{equation*}
where $i_1,\ldots,i_j\in\mathbb N\cup\{0\}$ and the constants $\hat\varphi_i$
are the Taylor series coefficients of $\log(z/\mathrm{sinh}\,z)$ 
as outlined in the proof.
\end{theorem}
\begin{proof}
We prove the result in three steps. First, we rewrite the 
density function $\phi$ in the form
\begin{equation*}
\phi(x)=\frac{1}{2\pi^2}\int_{\mathbb{R}}f(xz/\pi)\,
\exp\bigl(P\,\varphi(z)\bigr)\,\mathrm{d}z,
\end{equation*}
where for all $z\in\mathbb R$ we set 
$f(z)\coloneqq \cos(z)$ and 
$\varphi(z)\coloneqq\log\bigl(z/\mathrm{sinh}(z)\bigr)$.
We observe that $f$ and $\varphi$ are even functions and 
have power series expansions about $z=0$ (with infinite radii
of convergence) of the form 
\begin{equation*}
f(xz/\pi)=\sum_{k\geqslant0}\hat f_k\,z^{2k}
\qquad\text{and}\qquad 
\varphi(z)=-\tfrac16\,z^2+z^4\sum_{k\geqslant0}\hat\varphi_k\,z^{2k}
\end{equation*}
where $\hat f_k=(-1)^k(x/\pi)^{2k}/(2k)!$ and $\hat\varphi_0=1/180$, 
$\hat\varphi_1=-1/2835$ and so forth.
The coefficients $\hat\varphi_k$ can be analytically computed via
the Taylor coefficients of $\varphi$ to any order. In practice 
we computed them via Maple. We separate the 
quadratic term $-z^2/6$ from the series expansion for 
$\varphi(z)$ in $\exp\bigl(P\,\varphi(z)\bigr)$ and set
\begin{equation*}
g(z)\coloneqq\exp\Bigl(P\,z^4\sum_{k\geqslant0}\hat\varphi_kz^{2k}\Bigr).
\end{equation*}
Expanding $g=g(z)$ as a power series in $z$ we find
\begin{equation*}
g(z)=\sum_{\ell\geqslant0}\hat g_\ell\,z^{2\ell}
\end{equation*}
where explicitly we see that $g_0=1$, $g_1=0$ and for all $\ell\geqslant2$
we have
\begin{equation*}
\hat g_\ell=\sum_{j=1}^{\llcorner\ell/2\lrcorner}\frac{P^j}{j!}c_{j,\ell}
\end{equation*}
with the constants $c_{j,\ell}$ as stated in the Theorem.
Second, we observe $\varphi=\varphi(z)$ has a global maximum
at $z=0$ and apply the Laplace method outlined in \citet[pp.~272--4]{BO}.
Hence to within exponentially small errors, we shrink the range of integration
in the integral representation for $\phi$ above to an asymptotically small
interval strictly containing the origin $z=0$. We replace $f$ and $\varphi$
by their power series expansions about $z=0$ above. Then we extend
the range of integration to the whole real line. Thus, as $P\to\infty$, we obtain
\begin{equation*}
\phi(x)\sim\frac{1}{2\pi^2}\int_{\mathbb{R}} 
\mathrm{e}^{-\tfrac16 Pz^2}\cdot\sum_{k\geqslant0}(\hat f\star \hat g)_k\,z^{2k}
\,\mathrm{d}z,
\end{equation*}
where 
\begin{equation*}
(\hat f\star \hat g)_k=\sum_{\ell=0}^k\hat f_{k-\ell}\,\hat g_{\ell}.
\end{equation*}
Third, using the substitution $z=(3/P)^{1/2}\tau$ and the identity
\begin{equation*}
\int_{\mathbb{R}} \mathrm{e}^{-\frac12\tau^2}\tau^{2k}\,\mathrm{d}\tau
\equiv(2\pi)^{1/2}(2k-1)(2k-3)\cdots(1),
\end{equation*}
we see that
\begin{equation*}
\phi(x)\sim\biggl(\frac{3}{2\pi^3P}\biggr)^{1/2}
\sum_{k\geqslant0}\frac{a_k}{P^k}\cdot(\hat f\star \hat g)_k,
\end{equation*}
where the constants $a_k$ are defined in the statement of the Theorem. 
Hence we see that 
\begin{equation*}
\phi(x)\sim\biggl(\frac{3}{2\pi^3P}\biggr)^{1/2}
\sum_{k\geqslant0}\frac{a_k}{P^k}\sum_{\ell=0}^k\frac{(-1)^{k-\ell}}{\bigl(2(k-\ell)\bigr)!}
\Bigl(\frac{x}{\pi}\Bigr)^{2(k-\ell)}\cdot\hat g_\ell.
\end{equation*}
Using the explicit form for $\hat g_\ell$ given above, separating
out the cases $\ell=1,2$ in the last sum and integrating
with respect to $x$, generates the stated series representation
for $\Phi=\Phi(x)$.
\end{proof}

\newpage

\subsection*{(4) Chebychev coefficients}\label{sec:tables}
In Tables~1 to 4 we give the coefficients of the 
Chebychev polynomial approximations for the inverse distribution 
function $\Phi^{-1}=\Phi^{-1}(u)$ for the sum of $P$ Logistic 
random variables. We constructed the polynomials for 
$P=10^3,10^4,10^5,10^6$. In each case, using that $\Phi^{-1}$
is odd about $u=1/2$, we split the domain $[1/2,1-10^{-12}]$ 
into three regions:  the \emph{central} $[1/2,u_1]$; 
\emph{middle} $(u_1,u_2]$ and 
\emph{tail} $(u_2,1-10^{-12}]$ \emph{regions}. 
For the Chebychev coefficients $c_n$ we use the notation of \citet[Section~5.8]{PTVF},
namely the approximating Chebychev polynomial of degree $N$ for $z\in[-1,1]$ has the
form
\begin{equation*}
\tfrac12 c_0+c_1T_1(z)+c_2T_2(z)+\cdots+c_NT_N(z),
\end{equation*}
where $T_n(z)=\cos(n\,\mathrm{arccos}\,z)$ for $n=1,\ldots,N$ are the 
degree $N$ Chebychev polynomials. We also
quote the constants $k_1$ and $k_2$ used to ensure $z=-1$ and
$z=+1$ at the left and right endpoints, respectively, of the three regions 
(see Section~2 in the main manuscript for details). 
All coefficients were computed in Maple using $25$--$45$ digit accuracy and
then imported to Matlab for the Monte--Carlo simulations we performed (in double
precision arithmetic). Hence in the tables we quote the coefficients 
to double precision accuracy.

\begin{table}[ht]\label{p3table}
\caption{Case $P=10^3$: Chebychev coefficients $c_n$}
\begin{center}
\begin{tabular}{cccc}\hline\hline
$\phantom{\hat{\Big|}}$                                   
$n$ & central          &   middle                      &      tail     \\\hline
0   &  2.119420458542864e+00  &   1.477204569401002e+02  &   5.017891906926475e+02\\
1   &  6.366541597036217e-02  &   2.512841965456320e+01  &   1.527347358616282e+02\\
2   &  4.107853024715088e-03  &   1.130257188584296e+00  &   1.126378230488515e+00\\
3   &  3.294158206357919e-04  &  -1.282047839055625e-01  &  -4.705557379759551e-01\\
4   &  2.930679509853143e-05  &   9.575205785610623e-03  &   1.712745827810646e-01\\
5   &  2.770817602734147e-06 &   -2.256216011606160e-04  &  -5.254361417890216e-02\\
6   &  2.726390206722425e-07 &   -6.001369439751175e-05  &   1.460434532492525e-02\\
7   &  2.759420181388308e-08 &    1.075449360355294e-05  &  -3.718279846311429e-03\\
8   &  2.852054021544598e-09  &  -8.617168937931277e-07  &   8.606043517061642e-04\\
9   &  2.995958827063737e-10  &  -1.953300581656130e-09  &  -1.755629646754982e-04\\
10  &   3.187964787435130e-11 &    1.089315305113562e-08&     2.877465568099595e-05\\
11  &   3.428078463276479e-12 &   -1.603001402911367e-09 &   -2.372520758491911e-06\\
12  &   3.718011102620834e-13 &    1.004148628171250e-10 &   -7.636885527600424e-07\\
13  &   4.014002640668204e-14 &    5.955834555803564e-12 &    5.312592968146213e-07\\
14  &                         &   -2.253487007195484e-12 &   -2.073627141166733e-07\\
15  &                          &                         &    6.501618425527683e-08\\
16  &                          &                         &   -1.753681966723967e-08\\
17 &                           &                         &    4.081816622919563e-09\\
18 &                           &                         &   -7.756982638097922e-10\\
19 &                           &                         &    9.261438759264053e-11\\
20 &  $u_1$                    &    $u_2$                &    1.007783529912487e-11\\
21 & 8.083481113027166e-01     &  9.593726184247793e-01  &   -1.187272835871029e-11\\
22 &                           &                         &    5.249458333523391e-12\\
23 &                           &                         &   -1.771544847356913e-12\\
24 &                           &                         &    4.830855320498356e-13\\
   &   $k_1$                   &  $k_1$                  &    $k_1$                   \\
   &  1.017628007780257e-03    & 3.200822102223405e-02   &    6.587829001812997e-03   \\
   &   $k_2$                   &  $k_2$                  &    $k_2$                   \\
   &   -1.0                    & -2.614268030790329e+00  &   -1.743877010128950e+00   \\
\hline\hline
\end{tabular}
\end{center}
\end{table}

\begin{table}\label{table:p4}
\caption{Case $P=10^4$: Chebychev coefficients $c_n$}
\begin{center}
\begin{tabular}{cccc}\hline\hline
$\phantom{\hat{\Big|}}$                                   
$n$ & central                &   middle                &      tail     \\\hline
0   &  2.108319108843320e+00 &    4.423384392229279e+02 &    1.567580813837419e+03\\
1   &  5.728783093023983e-02 &    7.541649685957418e+01 &    4.900726813089462e+02\\
2   &  3.336830714052909e-03 &    3.802545600402460e+00 &    3.620608672934667e+00\\
3   &  2.414072718500087e-04 &   -4.162811951689960e-01 &   -1.753263346614506e+00\\
4   &  1.937066863959366e-05 &    2.847552228467595e-02 &    6.344679729235876e-01\\
5   &  1.651565758464117e-06 &   -2.504577026089741e-04 &   -1.987686632758167e-01\\
6   &  1.465393887537522e-07 &   -2.386566446362687e-04 &    5.635994736542883e-02\\
7   &  1.337341473177811e-08 &    3.491625902183290e-05 &   -1.459248555266122e-02\\
8   &  1.246313142218848e-09 &   -2.120845558475067e-06 &    3.413734029250331e-03\\
9   &  1.180432898053472e-10 &   -1.182684032497276e-07 &   -6.944358877927228e-04\\
10  &   1.132527541500540e-11&     4.322141511114120e-08 &    1.087996218632490e-04\\
11  &   1.098023528922786e-12&    -4.738539238127797e-09 &   -5.594596153804388e-06\\
12  &   1.073779033277571e-13&     1.260115332003234e-10 &   -4.896715048489821e-06\\
13  &   1.047559090342857e-14&     4.473491268259253e-11 &    2.847457352463447e-06\\
14  &                       &    -8.361451977660587e-12 &   -1.074062330827547e-06\\
15  &                       &                           &     3.323610314467442e-07\\
16  &                       &                          &   -8.847320546057122e-08\\
17  &                       &                          &    2.001510154092221e-08\\
18  &                       &                          &   -3.502074323316843e-09\\
19  &                       &                          &    2.603633442997274e-10\\
20  &                       &                          &    1.372771337705810e-10\\
21  &   $u_1$               &      $u_2$               &   -9.139416507192849e-11\\
22  &  7.958822967393328e-01&  9.509351131348488e-01   &    3.665796249029507e-11\\
23  &                       &                          &   -1.186014171700694e-11\\
24  &                       &                          &    3.279572949450883e-12\\
25  &                       &                          &   -7.546639202675646e-13\\
    &   $k_1$               &  $k_1$                   &    $k_1$                   \\
    & 1.105181672438157e-04 & 1.034459423699480e-02    &   2.045267156101968e-03  \\
    &   $k_2$               &  $k_2$                   &    $k_2$                   \\
    &   -1.0                &  -2.509336090190907e+00  &    -1.693843536107212e+00  \\
\hline\hline
\end{tabular}
\end{center}
\end{table}

\begin{table}\label{table:p5}
\caption{Case $P=10^5$: Chebychev coefficients $c_n$}
\begin{center}
\begin{tabular}{cccc}\hline\hline
$\phantom{\hat{\Big|}}$                                   
$n$ & central                 &   middle                &      tail     \\\hline
0   &  2.119671622170969e+00  &   1.477218788406007e+03 &    5.011286884624496e+03\\
1  &   6.363975712024890e-02  &   2.512862797352757e+02 &    1.522850796149265e+03\\
2   &  4.104766114250549e-03  &   1.129540906621810e+01 &    9.932064400088189e+00\\
3   &  3.290685366779347e-04  &  -1.284133303934886e+00 &   -4.841779474373805e+00\\
4   &  2.926766548631051e-05  &   9.580643879371174e-02 &    1.713947282173348e+00\\
5   &  2.766375465005536e-06  &  -2.254084712019299e-03 &   -5.256368571812046e-01\\
6   &  2.721311932661770e-07  &  -6.011921114110920e-04 &    1.460808124814458e-01\\
7   &  2.753581314991415e-08  &   1.076895266512730e-04 &   -3.718798654742950e-02\\
8   &  2.845309692030776e-09  &  -8.626049863871663e-06 &    8.606223447393348e-03\\
9   &  2.988139671124323e-10  &  -2.015930438797571e-08 &   -1.755399457767294e-03\\
10  &  3.178872129261131e-11  &   1.091637269032094e-07 &    2.876288620585080e-04\\
11  &  3.417484364751520e-12  &  -1.605856150405759e-08 &   -2.368489059323497e-05\\
12 &   3.706169816000400e-13  &   1.005307561650781e-09 &   -7.648223353997908e-06\\
13  &  4.048302208916000e-14 &    5.978528890984932e-11 &    5.315305491340833e-06\\
14  &  4.396573842960000e-15  &  -2.269982787480031e-11 &   -2.074161512888172e-06\\
15  &                         &   2.702274060697561e-12 &    6.502342260867294e-07\\
16  &                         &                         &   -1.753665035039117e-07\\
17  &                         &                         &    4.081237348842598e-08\\
18  &                         &                         &   -7.754302812987447e-09\\
19  &                         &                         &    9.252527846608224e-10\\
20  &                         &                         &    1.010210658045239e-10\\
21  &    $u_1$                &    $u_2$                &   -1.187795751600622e-10\\
22  &   8.083217460069005e-01 &  9.593729171835723e-01  &    5.250467255871130e-11\\
23  &                         &                         &   -1.775522861434831e-11\\
24 &                          &                         &    5.088026187684781e-12\\
25 &                          &                         &   -1.219201410524395e-12\\
   &   $k_1$                  &  $k_1$                  &    $k_1$                   \\
   & 1.017802054606037e-05    & 3.200351484212014e-03   &  6.587833651006969e-04      \\
   &   $k_2$                  &  $k_2$                  &    $k_2$                   \\
   &   -1.0                   & -2.613743480459985e+00  &    -1.743878946550866e+00  \\
\hline\hline
\end{tabular}
\end{center}
\end{table}

\begin{table}\label{p6table}
\caption{Case $P=10^6$: Chebychev coefficients $c_n$}
\begin{center}
\begin{tabular}{cccc}\hline\hline
$\phantom{\hat{\Big|}}$                                   
$n$ & central                 &   middle                &      tail     \\\hline
0  &   2.108344599409140e+00  &   4.423387698816169e+03 &    1.567373755547282e+04\\
1  &   5.728549519703990e-02  &   7.541655863088429e+02 &    4.899305223831529e+03\\
2  &   3.336575808689480e-03  &   3.802379058616604e+01 &    3.577332189281591e+01\\
3  &   2.413813273461038e-04  &  -4.163430999031781e+00 &   -1.757858358774398e+01\\
4  &   1.936802680113490e-05  &   2.847677160003316e-01 &    6.345091624110355e+00\\
5  &   1.651294878124777e-06  &  -2.503283183664147e-03 &   -1.987772090765017e+00\\
6  &   1.465114294664739e-07  &  -2.386933971327419e-03 &    5.636162519699199e-01\\
7  &   1.337051302211130e-08  &   3.492033773769139e-04 &   -1.459275570545880e-01\\
8  &   1.246010659010858e-09  &  -2.120974972399575e-05 &    3.413756240516940e-02\\
9  &   1.180116451284370e-10  &  -1.183066241781009e-06 &   -6.944288378230288e-03\\
10 &    1.132195529071883e-11 &    4.322932122559763e-07&     1.087946967867817e-03\\
11 &    1.097675364328587e-12 &   -4.739210876391624e-08 &   -5.592669103359117e-05\\
12 &    1.073514671625467e-13 &    1.259986653651446e-09 &   -4.897320182152358e-05\\
13 &    1.057521381174667e-14 &    4.473545677363332e-10 &    2.847620441432345e-05\\
14 &    1.038065498093333e-15 &   -8.342672530339772e-11 &   -1.074100043264248e-05\\
15 &                          &    6.400515074460786e-12 &    3.323680348310634e-06\\
16 &                          &                         &   -8.847395067541618e-07\\
17 &                          &                        &    2.001495695896450e-07\\
18 &                          &                         &   -3.501945056707199e-08\\
19 &                          &                         &    2.603085338007325e-09\\
20 &                          &                         &    1.372951378171891e-09\\
21 &                          &                         &   -9.139877345296721e-10\\
22 &  $u_1$                   &         $u_2$           &    3.665779362863886e-10\\
23 & 7.958796098523839e-01    &  9.509348932922126e-01  &   -1.185719114843435e-10\\
24 &                          &                         &    3.275872579877628e-11\\
25 &                          &                         &   -7.697011508395845e-12\\
26 &                          &                         &    1.442404227758017e-12\\
   &   $k_1$                  &  $k_1$                  &    $k_1$                 \\
   &  1.105201744870294e-06   & 1.034445867947374e-03   &   2.045266247283142e-04  \\
   &   $k_2$                  &  $k_2$                  &    $k_2$                  \\
   &   -1.0                   & -2.509285608336680e+00  &   -1.693842339092087e+00 \\
\hline\hline
\end{tabular}
\end{center}
\end{table}

\subsection*{(5) Finite representation}\label{sec:finite}
The probability density function $\phi$ of the sum of $P$ independent
identically distributed Logistic random variables is given by
the $P$-fold convolution of the density function for one
Logistic random variable. We therefore anticipate $\phi$ to
have a finite form. Indeed it has. We prove this here via residue
calculus using the form for the density function given as the 
inverse Fourier transform of the characteritic function in the form
\begin{equation*}
\phi(x)\coloneqq\frac{1}{2\pi}\int_{\mathbb{R}}
\biggl(\frac{\pi z}{\mathrm{sinh}\,\pi z}\biggr)^P
\mathrm{e}^{-\mathrm{i}xz}\,\mathrm{d}z.
\end{equation*}
\begin{theorem}
The even probability density function $\phi=\phi(x)$ has the finite
representation (for $x>0$)
\begin{equation*}
\phi(x)=-\mathrm{i}^{P+1}\mathrm{e}^{-x}\sum_{k=0}^{P-1}C_{P-1-k}
\sum_{\ell=0}^kA_\ell B_{k-\ell}\,x^\ell
\cdot\sum_{j=1}^{P+\ell-k}\frac{j!\,D_{P+\ell-k,j}}{(1-(-1)^P\mathrm{e}^{-x})^{j+1}}.
\end{equation*}
Here the constants $A_k$ and $B_k=B_k(P)$ are given by
\begin{equation*}
A_k=\frac{(-\mathrm{i})^k}{k!}
\qquad\text{and}\qquad 
B_k=\frac{P!}{(P-k)!k!}\cdot(-\mathrm{i})^{-k}.
\end{equation*}
The constants $C_k=C_k(P)$ are shifted Taylor series coefficients of 
$(\mathrm{sinh}\,\pi z)^{-P}$ about $z=0$ and the constants 
$D_{k,j}$ solve a linear system of equations 
(see how both sets of constants are outlined in the proof).
\end{theorem}
\begin{proof}
We prove the result in four steps. First, we choose a closed contour 
$\mathcal C$ in the complex $z$-plane given by the the interval $[-R,R]$ 
on the real axis and a semi-circular arc on the lower half complex plane
of radius $R$. Then integrating in the clockwise direction we see
for $x>0$ we have
\begin{equation*}
\phi(x)\coloneqq\lim_{R\to\infty}\frac{1}{2\pi}\int_{\mathcal C}
\biggl(\frac{\pi z}{\mathrm{sinh}\,\pi z}\biggr)^P
\mathrm{e}^{-\mathrm{i}xz}\,\mathrm{d}z.
\end{equation*}
Here we have used that in the limit $R\to\infty$ the contribution
to the contour integral from the semi-circular arc is vanishingly 
small for $x>0$. Note that the integrand has a removable singularity
at $z=0$ and poles at $z=\pm\mathrm{i}n$ for all $n\in\mathbb N$.
Hence by the Cauchy Residue Theorem for $x>0$ we have 
\begin{equation*}
\phi(x)=-\mathrm{i}\sum_{n\in\mathbb{N}}
\mathrm{residue}\bigl(\mathrm{e}^{-\mathrm{i}xz}
(\pi z/\mathrm{sinh}\,\pi z)^P\colon z=-\mathrm{i}n\bigr).
\end{equation*}
Second, our goal now is to compute the coefficient of the  
pole term in the Laurent series expansion of the integrand 
about each $z=-\mathrm{i}n$. We fix $n$ for the 
moment and set $\zeta\coloneqq z+\mathrm{i}n$. 
The integrand above has a pole at $\zeta=0$. We rewrite
the regular numerator terms in the integrand as follows
\begin{equation*}
\mathrm{e}^{-\mathrm{i}xz}=\mathrm{e}^{-nx}\sum_{k\geqslant0}A_kx^k\zeta^k,
\qquad\text{and}\qquad
(\pi z)^P=\pi^P (-\mathrm{i}n)^P\sum_{k=0}^PB_kn^{-k}\zeta^k,
\end{equation*}
where the coefficients $A_k=(-\mathrm{i})^k/k!$ 
and $B_k=P!(-\mathrm{i})^{-k}/(P-k)!k!$ as stated in the Theorem.
Using that 
$\mathrm{sinh}(\pi(\zeta-\mathrm{i}n))=(-1)^n\mathrm{sinh}\,\pi\zeta$ 
we rewrite the denominator term as a factor as follows
\begin{equation*}
(\mathrm{sinh}\,\pi z)^{-P}=\zeta^{-P}\pi^{-P}(-1)^{nP}\sum_{k\geqslant0}C_k\zeta^k
\end{equation*}
where the constants $C_k=C_k(P)$ are defined by this relation
and are shifted Taylor coefficients for 
$(\mathrm{sinh}\,\pi\zeta)^{-P}$ about $\zeta=0$.
Combining these three factors then, modulo the `constant' factor 
$\mathrm{e}^{-nx}\cdot(-\mathrm{i})^P\cdot(-1)^{nP}$
and taking into account the factor $\zeta^{-P}$ from the
last term, we are interested in determining the coefficient
of $\zeta^{P-1}$ in the product 
\begin{equation*}
\sum_{k\geqslant0}A_kx^k\zeta^k\cdot\sum_{k=0}^PB_kn^{-k}\zeta^k\cdot\sum_{k\geqslant0}C_k\zeta^k
=\sum_{k\geqslant0}\bigl(\{Ax\}\star\{B/n\}\star C\bigr)_k\zeta^k.
\end{equation*}
Here, as in Section~\ref{sec:largesum}, we denote
$(A\star B)_k=\sum_{\ell=0}^kA_{\ell}\,B_{k-\ell}$
and by $\{Ax\}$ and $\{B/n\}$ the series of coefficients 
$\{A_0,A_1x,A_2x^2,\ldots\}$ and $\{B_0,B_1n^{-1},B_2n^{-2},\ldots\}$, 
respectively. In other words we wish to determine 
\begin{equation*}
\bigl(\{Ax\}\star\{B/n\}\star C\bigr)_{P-1}
\equiv\sum_{k=0}^{P-1}C_{P-1-k}\sum_{\ell=0}^kA_\ell B_{k-\ell}\,n^{\ell-k}x^\ell.
\end{equation*}
Hence for $x>0$ we have
\begin{equation*}
\phi(x)=-\mathrm{i}\sum_{n\geqslant1}
\bigl((-1)^{P}\mathrm{e}^{-x}\bigr)^n(-\mathrm{i}n)^P
\sum_{k=0}^{P-1}C_{P-1-k}\sum_{\ell=0}^kA_\ell B_{k-\ell}\,n^{\ell-k}x^\ell.
\end{equation*}
Third, for convenience we define $Y=Y(x)$ by
$Y(x)\coloneqq(-1)^P\mathrm{e}^{-x}$. Then we can write $\phi$
in the form
\begin{equation*}
\phi(x)=(-\mathrm{i})^{P+1}\,Y
\sum_{k=0}^{P-1}C_{P-1-k}\sum_{\ell=0}^kA_\ell B_{k-\ell}\,x^\ell\cdot
\sum_{n\geqslant0}(n+1)^{P+\ell-k}\,Y^n,
\end{equation*}
where we shifted the summation variable $n$ by $1$ in the last step.
Fourth, our goal now is to rewrite the final sum over $n$ as
a finite sum. We need two identities to achieve this.
These are the straightforward identity for any integer 
$m\geqslant0$,
\begin{equation*}
(1-Y)^{-(m+1)}\equiv\frac{1}{m!}\sum_{n\geqslant0}(n+1)(n+2)\cdots(n+m)Y^n,
\end{equation*}
and the expansion for any real $\xi$,
\begin{equation*}
\xi(\xi+1)(\xi+2)\cdots(\xi+m)=\sum_{j=1}^{m+1}b_{m+1,j}\xi^j.
\end{equation*}
This relation defines the constants $b_{m+1,j}$ with $b_{m+1,m+1}=1$---more 
explicitly they are given by 
$b_{m+1,j}=\sum a_{i_1}a_{i_2}\cdots a_{i_{m+1-j}}$,
where the $a_i\in\{1,2,\ldots,m\}$ and the sum is over all 
choices of $a_{i_1}$ to $a_{i_{m+1-j}}$ with no two terms the same.
We now expand $\xi^m$ in the form 
\begin{align*}
\xi^m=&\;D_{m,m}\,\xi(\xi+1)(\xi+2)\cdots(\xi+m)\\
&\;+D_{m,m-1}\,\xi(\xi+1)(\xi+2)\cdots(\xi+m-1)\\
&\;~~~~\vdots\\
&\;+D_{m,1}\,\xi.
\end{align*}
Observe that the constants $D_{m,1},\ldots,D_{m,m}$ solve the linear system of equations
\begin{equation*}
\begin{pmatrix}
1 & b_{2,1} & b_{3,1} & b_{4,1} & \cdots & b_{m,1} \\
0 & 1      & b_{3,2} & b_{4,2} & \cdots & b_{m,2} \\
0 & 0      & 1       & b_{4,3} & \cdots & b_{m,3} \\
\vdots&\vdots&\ddots & \ddots & \ddots & \vdots \\
0     &  0   & 0     &  0     &  1     & b_{m,m-1} \\
0     &  0   & 0     &  0     &  0     &  1  
\end{pmatrix}
\begin{pmatrix}
D_{m,1} \\ D_{m,2} \\ D_{m,3} \\ \vdots \\ D_{m,m-1} \\ D_{m,m}
\end{pmatrix}
=
\begin{pmatrix}
0 \\ 0 \\ 0 \\ \vdots \\ 0 \\ 1
\end{pmatrix}.
\end{equation*}
If we now take $\xi=n+1$ and using our identity for $(1-Y)^{-(m+1)}$ above 
we see 
\begin{align*}
\sum_{n\geqslant0}(n+1)^m\,Y^n
=&\;D_{m,m}\sum_{n\geqslant0}(n+1)(n+2)\cdots(n+m)\,Y^n\\
&\;+D_{m,m-1}\sum_{n\geqslant0}(n+1)(n+2)\cdots(n+m-1)\,Y^n\\
&~~~~\vdots\\
&\;+D_1\sum_{n\geqslant0}(n+1)\,Y^n\\
=&\;\sum_{j=1}^mj!\,D_{m,j}(1-Y)^{-(j+1)}.
\end{align*}
Taking $m=P+\ell-k$ and substituting back for 
$Y=(-1)^P\mathrm{e}^{-x}$ completes the proof.
\end{proof}







\begin{thebibliography}{00}

\bibitem[\protect\citeauthoryear{Bender and Orszag}{1999}]{BO}
Bender, C.M. and Orszag, S.A.,
\textit{Advanced mathemtical methods for scientists 
and engineers I: Asymptotic methods and perturbation theory.}
Springer, 1999.

\bibitem[\protect\citeauthoryear{Gaines and Lyons}{1994}]{GL1994} 
Gaines, J. and Lyons, T.J.,
\textit{Random generation of stochastic area integrals.}
SIAM J. Appl. Math.
\textbf{54}(4) (1994), 1132--1146.

\bibitem[\protect\citeauthoryear{Gallier and Xu}{2002}]{GX} 
Gallier, J. and Xu, D.,
\textit{Computing exponentials of skew symmetric matrices and logarithms
of orthogonal matrices.} Int. J. Robot. Autom. \textbf{17} (2002), 1-–11.

\bibitem[\protect\citeauthoryear{H\"ormann}{1993}]{H}
H\"ormann, W.,
\textit{The transformed rejection method for generating Poisson random variables.}
Insurance: Mathematics and Economics \textbf{12} (1993), 39--45. 

\bibitem[\protect\citeauthoryear{Kloeden \textit{et al.}}{1992}]{KPW} 
Kloeden, P.E., Platen, E. and Wright, W.,
\textit{The approximation of multiple stochastic integrals.}
Stochastic Anal. Appl. \textbf{10} (1992), 431--441.

\bibitem[\protect\citeauthoryear{L\'evy}{1951}]{Levy} 
L\'evy, P.,
\textit{Wiener's random function and other Laplacian random functions.} 
Second Symposium of Berkeley. Probability and Statistics, UC Press, (1951), 171--186.

\bibitem[\protect\citeauthoryear{Malham and Wiese}{2013}]{MW2013} 
Malham, S.J.A. and Wiese, A.,
\textit{Chi-square simulation of the CIR process and the Heston model.}
IJTAF \textbf{16}(3) (2013), 1--38.

\bibitem[\protect\citeauthoryear{Mansuy and Yor}{2008}]{MY}
Mansuy, R. and Yor, M.,
\textit{Aspects of Brownian motion.} Springer--Verlag, 2008.

\bibitem[\protect\citeauthoryear{Moro}{1995}]{Moro} 
Moro, B., \textit{The full Monte.}
Risk \textbf{8}(2) (1995), 57--58.

\bibitem[\protect\citeauthoryear{Press \textit{et al.}}{1992}]{PTVF} 
Press, W.H., Teukolsky, S.A., Vetterling, W.T.
and Flannery, B.P., \textit{Numerical recipes in C:
The art of scientific computing.} Second Edition,
Cambridge University Press, 1992.

\bibitem[\protect\citeauthoryear{Ryd\'en and Wiktorsson}{2001}]{RW} 
Ryd\'en, T. and Wiktorsson, M.,
\textit{On the simulation of iterated It\^o integrals.}
Stochastic Process. Appl. \textbf{91} (2001), 151--168.

\bibitem[\protect\citeauthoryear{Stump and Hill}{2005}]{SH} 
Stump, D.~M. and Hill, J.~M.,
\textit{On an infinite integral arising in the numerical
integration of stochastic differential equations.}
Proc. R. Soc. A \textbf{461} (2005), 397--413.

\bibitem[\protect\citeauthoryear{Wiktorsson}{2001}]{W} 
Wiktorsson, M.,
\textit{Joint characteristic function and simultaneous
simulation of iterated It\^o integrals for multiple independent Brownian
motions.} Ann. Appl. Probab. \textbf{11}(2) (2001), 470--487.

\end{thebibliography}
\end{document}